\documentclass{article} 
\usepackage[british]{babel}
\usepackage{makecell} 
\usepackage{array}
\emergencystretch=\maxdimen
\usepackage{indentfirst}
\usepackage[utf8]{inputenc}
\usepackage[T1]{fontenc}

\usepackage{hyperref} 
\usepackage[top=1in, bottom=1in, left=1in, right=1in]{geometry}

\usepackage{lipsum}

\usepackage{graphicx}
\usepackage{epstopdf}
\usepackage{epsfig}
\usepackage{amsthm} 
\usepackage{amsmath} %

\usepackage{amssymb}
\usepackage{mathtools} 
\usepackage[adobe-utopia]{mathdesign}

\usepackage{multirow}

\newenvironment{myproof}[1][\proofname]{%
  \begin{proof}[#1]$ $\par\nobreak\ignorespaces   
}{%

  \end{proof}
}

\theoremstyle{definition}
\newtheorem{definition}{Definition}[section] %

\newtheorem{example}{Example}[section]
\newtheorem{conclusion}{Conclusion}[section]

\newtheorem{theorem}{Theorem}[section]

\newcommand*{\titleGM}{\begingroup 
\begin{center}
\hbox{ 
\hspace*{0.05\textwidth} 
\hspace*{0.05\textwidth} 
\parbox[b]{0.98\textwidth}
{ 
{\Huge\bfseries Algebra of paravectors \\[0.5\baselineskip] }
\\[2\baselineskip] 
{\Large \textsc{Józef Radomański}} 
\\
ul. Bohaterów Warszawy 9, 48-300 Nysa, Poland
\\
e-mail: radomanski@wikakwa.pl

\vspace{0.05\textheight} 
}}

\end{center}
\begin{abstract}

\textit{Paravectors just like integers have a ring structure. By introducing an integrated product we get geometric properties which make paravectors similar to vectors. The concepts of parallelism, perpendicularity and the angle are conceptually similar to vector counterparts, known from the Euclidean geometry. Paravectors meet the idea of parallelogram law, Pythagorean theorem and many other properties well-known to everyone.}
\end{abstract}

Keywords: \textit{Paravector, multivector, Geometric Algebra}

\endgroup}

\begin{document}

\pagestyle{empty} 

\titleGM 

\pagestyle{plain} 


This article refers to professor William Baylis' researches and shows a surprising similarity between properties of paravectors and vectors in the Euclidean space. It can be confusing at first to get used to a column notation of paravectors, but this form has proved to be the most adequate. The main advantage of this notation is a transparency of mathematical transformations. The author hopes that in future paravector formalism will replace the multivectors formalism in physics as simpler, more intuitive and imaginable. The results presented in this article justify it fully.
\section{Basic definitions}
\begin{definition}
The term \textbf{paravector} means a pair consisting of a complex number ($\alpha$) and a vector ($\pmb{\beta}$) belonging to a three-dimensional complex space.
\begin{equation} \label{eq:3.1}
\Gamma \eqqcolon
\begin{bmatrix}
\alpha \\ 
\pmb{\beta }
\end{bmatrix}
=\begin{bmatrix}
a+id \\ 
\mathbf{b}+i\mathbf{c}
\end{bmatrix}
\end{equation}
\end{definition}
The number will be called a \textbf{scalar}. Paravectors will be denoted with capital letters, eg: $ A, X, \Gamma $. Greek letters will mean a complex size, and Roman letters - a real one. A scalar component of paravector $\Gamma$ will be denoted with index ,,S'', and a vector component with ,,V'', i.e.: $\Gamma_S =\alpha$ and $\Gamma_V =\pmb{\beta}$.
\begin{definition}
The \textbf{reversed} element of paravector (1) is the paravector
\begin{equation}
\Gamma ^{-}
\eqqcolon
\begin{bmatrix}
a+id \\ 
-\mathbf{b-}i\mathbf{c}
\end{bmatrix}
\end{equation}
\end{definition}
\begin{definition}
The \textbf{conjugated} element of paravector (1) is the paravector 
\begin{equation}
\Gamma ^{\ast}
\eqqcolon
\begin{bmatrix}
a-id \\ 
\mathbf{b-}i\mathbf{c}
\end{bmatrix}
\end{equation}
\end{definition}
\begin{definition}
of the \textbf{summation}:
\begin{equation}
\begin{bmatrix}
\alpha _{1} \\ 
\pmb{\beta }_{1}
\end{bmatrix}
+\begin{bmatrix}
\alpha _{2} \\ 
\pmb{\beta }_{2}
\end{bmatrix}
\eqqcolon 
\begin{bmatrix}
\alpha _{1}+\alpha _{2} \\ 
\pmb{\beta }_{1}+\pmb{\beta }_{2}
\end{bmatrix}
\end{equation}
\end{definition}
\begin{conclusion}The \textbf{neutral} element under addition (\textbf{null element}) is the paravector
\begin{equation}
0\eqqcolon 
\begin{bmatrix}
0 \\ 
\mathbf{0}
\end{bmatrix}
\end{equation}
\end{conclusion}
\begin{definition}
The \textbf{opposite} element of paravector (1) with respect to the addition is 
\begin{equation}
-\Gamma
\eqqcolon
\begin{bmatrix}
-\alpha \\ 
-\pmb{\beta }
\end{bmatrix}
\end{equation}
\end{definition}
\begin{definition} of the \textbf{multiplication}: \label{def:1.6}
\begin{equation}
\begin{bmatrix}
\alpha _{1} \\ 
\pmb{\beta }_{1}
\end{bmatrix}
\begin{bmatrix}
\alpha _{2} \\ 
\pmb{\beta }_{2}
\end{bmatrix}%
\eqqcolon 
\begin{bmatrix}
\alpha _{1}\alpha _{2}+\pmb{\beta }_{1}\pmb{\beta }_{2} \\ 
\alpha _{2}\pmb{\beta }_{1}+\alpha _{1}\pmb{\beta }_{2}+i\pmb{\beta 
}_{1}\times \pmb{\beta }_{2}
\end{bmatrix}
\end{equation}
where $\pmb{\beta }_{1}\pmb{\beta }_{2}$ is a scalar (dot) product of vectors, and $\pmb{\beta }_{1}\times \pmb{\beta }_{2}$ is a vector (cross) product.
\end{definition}
\begin{conclusion}
The \textbf{neutral} element under multiplication is the paravector
\begin{equation}
1 \eqqcolon 
\begin{bmatrix}
1 \\ 
\mathbf{0}
\end{bmatrix}
\end{equation}
\end{conclusion}
\noindent \textbf{Note:} There is no difference if we write number $\alpha $ or paravector 
$\begin{bmatrix}
\alpha \\ 
\mathbf{0}
\end{bmatrix}$.
\begin{conclusion}
The operation of multiplication is associative but not commutative, because the vector product is non-commutative.
\end{conclusion}
\begin{definition}
We call an outcome of multiplication of any paravector $\Gamma$ by the element conjugate $\Gamma^{\ast}$ the \textbf{vigor of paravector}:
\begin{equation}
\text{vig}\Gamma 
 \eqqcolon
\Gamma \Gamma ^{\ast }
\end{equation}
\end{definition}
\begin{conclusion}
To each paravector we can assign a vigor which is a real paravector and its scalar component is a positive number.
\end{conclusion}
\begin{proof}
$$
\vspace{1pt}
\Gamma \Gamma^{\ast}=
\begin{bmatrix}
\alpha \\ 
\pmb{\beta }
\end{bmatrix}
\begin{bmatrix}
\alpha^{\ast} \\ 
\pmb{\beta^{\ast} }
\end{bmatrix}
= 
\begin{bmatrix}
\alpha\alpha^{\ast}+\pmb{\beta}\pmb{\beta}^{\ast}\\ 
\alpha\pmb{\beta }^{\ast}+\alpha^{\ast}\pmb{\beta}
+i\pmb{\beta } \times \pmb{\beta }^{\ast}
\end{bmatrix}=
\begin{bmatrix}
(a +id)(a -id) +
 (\mathbf{b}+i\mathbf{c})(\mathbf{b}-i\mathbf{c})\\ 
(a +id)(\mathbf{b}-i\mathbf{c})+(a -id)(\mathbf{b}+i\mathbf{c})
+ i(\mathbf{b}+i\mathbf{c })\times (\mathbf{b}-i\mathbf{c })
\end{bmatrix}=
$$
$$
=\begin{bmatrix}
a^{2}+b^{2}+c^{2}+d^{2} \\ 
2(a\mathbf{b} + d\mathbf{c } + \mathbf{b} \times \mathbf{c })
\end{bmatrix}
$$
\end{proof}
\begin{definition}
We call an outcome of multiplication of any paravector $\Gamma$ by the reverse element $\Gamma^{-}$ the \textbf{determinant of a paravector}
\begin{equation}
\ \det \Gamma \eqqcolon \Gamma \Gamma
^{-}=\Gamma ^{-}\Gamma 
\end{equation}
\end{definition}
\begin{conclusion}
Each paravector has a determinant which is a complex number.
\end{conclusion}
\begin{proof}
$$
\vspace{1pt}
\Gamma \Gamma^{-}=
\begin{bmatrix}
\alpha \\ 
\pmb{\beta }
\end{bmatrix}
\begin{bmatrix}
\alpha \\ 
-\pmb{\beta }
\end{bmatrix}
= 
\begin{bmatrix}
\alpha^{2} - \beta^{2}\\ 
\alpha\pmb{\beta }-\alpha\pmb{\beta }
-i\pmb{\beta } \times \pmb{\beta }
\end{bmatrix}=
\begin{bmatrix}
(a +id)^{2} -
 (\mathbf{b}+i\mathbf{c})^{2}\\ 
-i(\mathbf{b}+i\mathbf{c }) \times (\mathbf{b}+i\mathbf{c })
\end{bmatrix}=
\begin{bmatrix}
a^{2}-b^{2}+c^{2}-d^{2}+2i(ad-\mathbf{b}\mathbf{c }) \\ 
0
\end{bmatrix}
$$
\end{proof}
\begin{conclusion}
The reversion and conjugation have the following properties:
  
  \begin{tabular}{|c|c|c|}
        \hline  
         &  Reversion of paravector & Conjugation of paravector \\
         \hline	\hline
         1 & $\left( \Gamma^{-}\right) ^{-}=\Gamma$ & $\left( \Gamma^{\ast }\right) ^{\ast }=\Gamma$ \\
         \hline
           2 & $\left( \Gamma+\Psi\right) ^{-}=\Gamma^{-}+\Psi^{-}$ & $\left( \Gamma+\Psi\right) ^{\ast }=\Gamma^{\ast}+\Psi^{\ast}$ \\
         \hline
           3  & $\left( \Gamma \Psi\right) ^{-}=\Psi^{-}\Gamma^{-}$ & $\left( \Gamma\Psi\right) ^{\ast }=\Psi^{\ast}\Gamma^{\ast }$ \\
         \hline
           4 & $\Gamma \Gamma^{-} \in C$ &  $\Gamma \Gamma^{\ast }\in R_{+}\times R^{3}$\\
           \hline
           5 & \multicolumn{2}{|c|}{$\left( \Gamma^{-} \right)^{\ast} = \left( \Gamma^{\ast} \right)^{-}$} \\
           \hline        
        \end{tabular}
\end{conclusion}
\begin{conclusion}
A set of pravectors together with an operation of summation forms an Abelian group, and with multiplication forms a semigroup. Therefore, we can conclude that a set of paravectors together with operations of summation and multiplication gives a ring with multiplicative identity.
\end{conclusion}
\begin{definition} \label{def:2.1.9}
We call the paravector $\Gamma $ \textbf{proper} if $\det \Gamma \in R_{+}\setminus\lbrace0\rbrace$ (the determinant is a positive real number).
\end{definition}
\begin{definition}
We call the paravector $\Gamma $ \textbf{singular} if $\det \Gamma =0$.
\end{definition}
By definition of the determinant it follows:
\begin{conclusion}
Each proper or singular paravector \eqref{eq:3.1}  must fulfill the following condition:
\begin{equation}
ad=\mathbf{bc}
\end{equation}
\end{conclusion}
\noindent \textbf{It is advisable to remember the above equation}, because it will be used with many proofs on proper paravectors.
\begin{conclusion}
Let $\Gamma _{1}$, $\Gamma _{2}$ be paravectors, then the following statements are true:
\begin{itemize}
\item $\det (\Gamma _{1}\Gamma _{2})=\det \Gamma _{1}\cdot \det \Gamma _{2}$
\item $\det (\Gamma ^{-})=\det \Gamma $
\item $\det (\Gamma ^{\ast })=(\det \Gamma)^{\ast} $
\end{itemize}
\end{conclusion}
\begin{definition}
For each non-singular paravector $\Gamma $, there exists an \textbf{inverse} element under multiplication:
\begin{equation}
\Gamma ^{-1} \eqqcolon \frac{\Gamma ^{-}}{\det \Gamma }
\end{equation}
\end{definition}
\begin{conclusion}
A set of non-singular paravectors together with multiplication is a non-commutative group.
\end{conclusion}
\begin{conclusion}A set of proper paravectors together with multiplication is a non-commutative group.
\end{conclusion}
\begin{definition} \label{def:2.1.12}
For each proper and/or singular paravector we define the \textbf{module} of paravector:
\begin{equation}
\left\vert \Gamma \right\vert  \eqqcolon \sqrt{\det \Gamma }
\end{equation}
\end{definition}
\begin{conclusion}
The module of paravector (proper or singular!) satisfies the following conditions:
\begin{enumerate}
\item $\left\vert s\Gamma \right\vert =\left\vert s\right\vert \left\vert
\Gamma \right\vert \qquad $ where $s\in R$
\item $\left\vert \Gamma _{1}\right\vert \left\vert \Gamma _{2}\right\vert
=\left\vert \Gamma _{1}\Gamma _{2}\right\vert $
\end{enumerate}
\end{conclusion}
\begin{definition} \label{def:2.1.13}
We call the paravector $\Lambda $ \textbf{orthogonal} if $\det \Lambda =1$ , or equivalently:
\begin{equation}
\Lambda  \eqqcolon \frac{\Gamma }{\sqrt{\det \Gamma }},\qquad \text{where } \Gamma \text{ is a proper paravector}
\end{equation}
\end{definition}
Directly from the above definition it follows:
\begin{conclusion}
If $\Lambda $ is an orthogonal paravector, then 
$\Lambda ^{-1}=\Lambda ^{-}$.
\end{conclusion}
\begin{definition}
We call the paravector $\Gamma$ \textbf{special} if $\Gamma^{-}=\Gamma^{\ast}$, or equivalently:
\begin{equation}
\nonumber
\begin{bmatrix}
a \\ 
i\mathbf{c}
\end{bmatrix}
,\qquad \text{where } a\in R, \text{ and } \mathbf{c}\in R^{3}
\end{equation}
\end{definition}
\begin{definition}
We call the paravector $\Gamma$ \textbf{unitar} if  $\Gamma\Gamma^{\ast}=1$
\end{definition}
To summarize the current knowledge about paravectors, we can say that very little is missing so that a~set of paravectors with summation and multiplication operations is a field: multiplication of paravectors is not commutative, and the role of the null element under multiplication is played by singular paravectors. Multiplying any paravector by singular one, we get a singular paravector. Note that although there are many null elements under multiplication, only one element is neutral with respect to summation.
\begin{conclusion}
The set of special paravectors together with operations summation and multiplication is a~division ring.
\end{conclusion}
\vspace{1pt}In the diagram (fig. 1) a ring of paravectors with some of its substructures is presented:
\begin{figure}[h]
\centering
\includegraphics[width=1\textwidth]{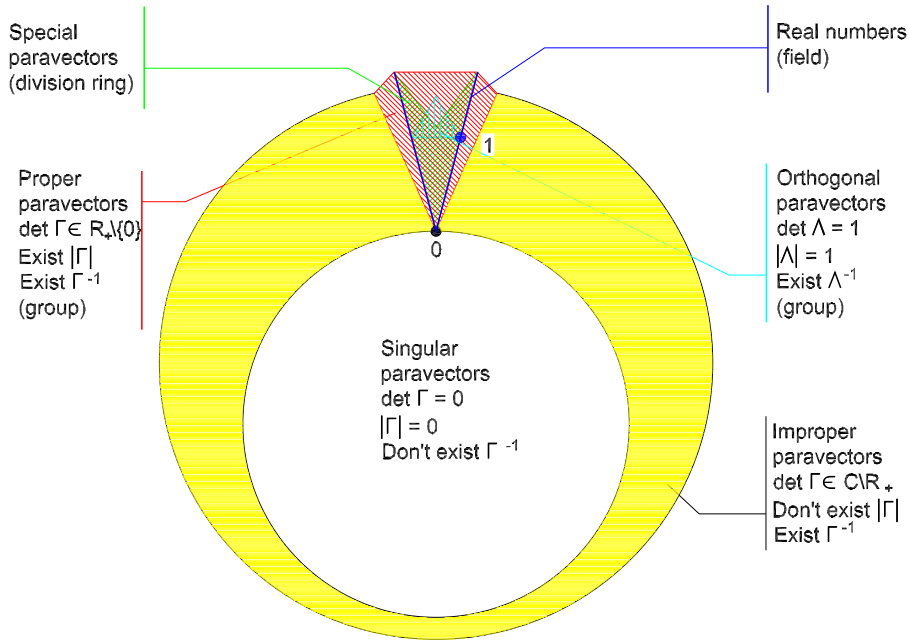}
\caption{The ring of paravectors with some substructures}
\end{figure}

\vspace{1pt}- A set of proper paravectors together with multiplication is a noncommutative group.

\vspace{1pt}- A set of orthogonal paravectors together with multiplication is a subgroup of proper paravectors group.

\vspace{1pt}The black point means zero, and the blue point is one.

\vspace{1pt}To avoid misunderstandings, we recall the following names which we will use:
\begin{itemize}
\item \vspace{1pt}$\Gamma ^{-}$ - The paravector reverse of $\Gamma$
\item \vspace{1pt}$\Gamma ^{-1}$ - The paravector inverse of $\Gamma$
\end{itemize}

\section{Integrated product of paravectors}
Synge \cite{Synge} defines the scalar product of complex quaternions as $\left( \Gamma_{1}\Gamma _{2}^{-}+\Gamma _{2}\Gamma _{1}^{-}\right) /2$, Hestenes \cite{Hestenes} as a scalar part of the product of multivectors $\Gamma _{1}\Gamma _{2}$. Analyzing the expression $\Gamma _{1} \Gamma _{2}^{-}$, we can see that it plays a similar but more universal role in the paravectors algebra as the dot product in vector space. Properties of its scalar component are the same as properties of the scalar product of vectors, and its vector part as a cross product of vectors. The trouble is that there can be two different products which have the same properties (the second one is $\Gamma _{1}^{-}\Gamma _{2}$), so we define two integrated products:

\begin{definition}
For any paravectors $\Gamma _{1}$ and $\Gamma _{2}$, the \textbf{right integrated product} of paravectors $\Gamma _{1}$ and $\Gamma _{2}$, denoted $\left( \Gamma _{1},\Gamma_{2}\right\rangle$, is defined by
\begin{equation} \nonumber
\left( \Gamma _{1},\Gamma
_{2}\right\rangle  \eqqcolon \Gamma _{1}\Gamma _{2}^{-}
\end{equation}
\end{definition}

hence $\qquad \left( \Gamma _{1},\Gamma _{2}\right\rangle = \Gamma _{1}\Gamma _{2}^{-}=
            \begin{bmatrix}
				\alpha _{1} \\ 
				\pmb{\beta }_{1}
				\end{bmatrix}
				\begin{bmatrix}
                \alpha _{2} \\ 
                -\pmb{\beta }_{2}
                \end{bmatrix}
                = 
                \begin{bmatrix}
                \alpha _{1}\alpha _{2}-\pmb{\beta }_{1}\pmb{\beta }_{2} \\ 
                -\alpha _{1}\pmb{\beta }_{2}+\alpha _{2}\pmb{\beta }_{1}-i\pmb{\beta 
                }_{1}\times \pmb{\beta }_{2}
                \end{bmatrix}$

\begin{definition}
For any paravectors $\Gamma _{1}$ and $\Gamma _{2}$ the \textbf{left integrated product} of paravectors $\Gamma _{1}$ and $\Gamma _{2}$, denoted $\left\langle \Gamma
_{1},\Gamma _{2}\right)$, is defined by
\begin{equation} \nonumber
\left\langle \Gamma
_{1},\Gamma _{2}\right)  \eqqcolon \Gamma _{1}^{-}\Gamma _{2}
\end{equation}
\end{definition}

or $\qquad \left\langle \Gamma _{1},\Gamma _{2}\right) = \Gamma _{1}^{-}\Gamma _{2}=
				\begin{bmatrix}
                \alpha _{1} \\ 
				-\pmb{\beta }_{1}
				\end{bmatrix}
				\begin{bmatrix}
                \alpha _{2} \\ 
                \pmb{\beta }_{2}
                \end{bmatrix}
                = 
                \begin{bmatrix}
                \alpha _{1}\alpha _{2}-\pmb{\beta }_{1}\pmb{\beta }_{2} \\ 
                \alpha _{1}\pmb{\beta }_{2}-\alpha _{2}\pmb{\beta }_{1}-i\pmb{\beta 
                }_{1}\times \pmb{\beta }_{2}
                \end{bmatrix}$
                
\noindent In both cases, the scalar part is the same, then
\begin{definition}\label{def:2.3}
The scalar component of an integrated product of paravectors 
$\Gamma _{1}$ and $\Gamma _{2}$, is called a \textbf{scalar product} of paravectors.
The scalar product will be denoted
\begin{equation}
\nonumber \left\langle \Gamma _{1},\Gamma
_{2}\right\rangle \eqqcolon ( \Gamma _{1},\Gamma_{2}\rangle_{S}=
\langle \Gamma _{1},\Gamma_{2} )_{S}.
\end{equation}
\end{definition}
\begin{definition}
The vector component of an integrated product we call a 
\textbf{vector product} of paravectors.
\begin{equation}
\nonumber ( \Gamma _{1},\Gamma_{2}\} \eqqcolon ( \Gamma _{1},\Gamma_{2}\rangle_{V}= -\langle \Gamma _{1},\Gamma_{2} )_{V}^{\ast}
\end{equation} 
\end{definition}
It is necessary to distinguish the orientation of a vector product, the same as at integrated product case. Therefore, we denote the right vector product
$( \Gamma _{1},\Gamma_{2}\}$, and the left vector product: $\{ \Gamma _{1},\Gamma_{2})$

\vspace{1pt}For further consideration it does not matter much if a product is right or left, so talking about an integrated product we will mean the right product. As well, this may be the left product, but it is important to use only one product constantly.
Hence the integrated product can be denoted as follow:
\begin{equation}
\left( \Gamma_{1}, \Gamma_{2} \right\rangle = \begin{bmatrix}
\left\langle \Gamma_{1}, \Gamma_{2} \right\rangle \\
( \Gamma_{1}, \Gamma_{2} \}
\end{bmatrix}=
\begin{bmatrix}
\left( \Gamma_{1}, \Gamma_{2} \right\rangle_{S} \\
\left( \Gamma_{1}, \Gamma_{2} \right\rangle_{V}
\end{bmatrix}
\end{equation}
and
\begin{equation}\label{eq:2.16}
\text{det} (\Gamma_{1}, \Gamma_{2} \rangle = \left\langle \Gamma_{1}, \Gamma_{2} \right\rangle^{2} - ( \Gamma_{1}, \Gamma_{2} \}^{2} =
\left( \Gamma_{1}, \Gamma_{2} \right\rangle^{2}_{S} - \left( \Gamma_{1}, \Gamma_{2} \right\rangle^{2}_{V} = \text{det} \Gamma_{1} \text{det} \Gamma_{2}
\end{equation}

\begin{theorem}
An integrated product of paravectors has the following properties:

\begin{tabular}{|l|l|l|}
\hline
& Right integrated product & Left integrated product \\ \hline\hline
1 & \multicolumn{2}{|c|}{Integrated product is a paravector} \\ \hline
2 & $\left( \Gamma _{1}+\Gamma _{2},\Gamma _{3}\right\rangle =\left( \Gamma
_{1},\Gamma _{3}\right\rangle +\left( \Gamma _{2},\Gamma _{3}\right\rangle $
& $\left\langle \Gamma _{1}+\Gamma _{2},\Gamma _{3}\right) =\left\langle
\Gamma _{1},\Gamma _{3}\right) +\left\langle \Gamma _{2},\Gamma _{3}\right) $
\\ \hline
3 & $\left( \alpha \Gamma _{1},\Gamma _{2}\right\rangle =\alpha \left(
\Gamma _{1},\Gamma _{2}\right\rangle =\left( \Gamma _{1},\alpha \Gamma
_{2}\right\rangle $ & $\left\langle \alpha \Gamma _{1},\Gamma _{2}\right)
=\alpha \left\langle \Gamma _{1},\Gamma _{2}\right) =\left\langle \Gamma_{1},\alpha \Gamma _{2}\right) $ \\ \hline
4 & $\left( \Gamma _{1},\Gamma _{2}\right\rangle ^{-}=\left( \Gamma
_{2},\Gamma _{1}\right\rangle $ & $\left\langle \Gamma _{1},\Gamma
_{2}\right) ^{-}=\left\langle \Gamma _{2},\Gamma _{1}\right) $ \\ \hline
5 & \multicolumn{2}{|l|}{$\left( \Gamma ,\Gamma \right\rangle =\left\langle
\Gamma ,\Gamma \right) =\det \Gamma \in C$} \\ \hline
6 & \multicolumn{2}{|l|}{$\det \left( \Gamma _{1},\Gamma _{2}\right\rangle
=\det \left\langle \Gamma _{1},\Gamma _{2}\right) =\det \Gamma _{1}\det
\Gamma _{2}$} \\ \hline
\end{tabular}
\end{theorem}

\begin{conclusion}
Let $\Gamma _{1},\Gamma _{2}$ and $\Gamma _{3}$ be any paravectors, then the following table shows the properties of the scalar product of paravectors compared to the properties of the scalar product of vectors in Euclidean space:

\begin{tabular}{|l|l|l|}
\hline 
& Scalar product of paravectors & Scalar product of vectors \\ \hline\hline
1 & $\left\langle \Gamma _{1},\Gamma _{2}\right\rangle \in C$ & $\left\langle \mathbf{x}_{1},\mathbf{x}_{2}\right\rangle \in R$ \\ \hline
2 & $\left\langle \Gamma _{1}+\Gamma _{2},\Gamma _{3}\right\rangle
=\left\langle \Gamma _{1},\Gamma _{3}\right\rangle +\left\langle \Gamma
_{2},\Gamma _{3}\right\rangle $ & same \\ \hline
3 & $\left\langle \alpha \Gamma _{1},\Gamma _{2}\right\rangle =\alpha
\left\langle \Gamma _{1},\Gamma _{2}\right\rangle =\left\langle \Gamma
_{1},\alpha \Gamma _{2}\right\rangle $ & same \\ \hline
4 & $\left\langle \Gamma _{1},\Gamma _{2}\right\rangle =\left\langle \Gamma
_{2},\Gamma _{1}\right\rangle $ & same \\ \hline
5 & $\left\langle \Gamma ,\Gamma \right\rangle \in C$ & $\left\langle 
\mathbf{x},\mathbf{x}\right\rangle \in R_{+}$ \\ \hline
6 & If $\left\langle \Gamma ,\Gamma \right\rangle =0$, then $\Gamma $
is singular & If $\left\langle \mathbf{x},\mathbf{x}\right\rangle
=0$, then $\mathbf{x}=0$ \\ \hline
\end{tabular}
\end{conclusion}

Rows 5 and 6 of Table 2.1 show a fundamental difference between the definition of the scalar product of paravectors and known algebraic definitions. These properties make me wonder whether to give another name for scalar product. However, the traditional name is left, because:
\begin{enumerate}
\item the value of the product is a scalar,
\item the role which the scalar product plays in the Paravector Algebra corresponds completely to the role of the dot product of vectors in Euclidean Geometry,
\item for spatial paravectors ($\Gamma_S = 0$), the scalar product of paravectors becomes the scalar product of vectors.
\end{enumerate}

\section{Geometrical properties of paravectors}
\vspace{1pt}By checking up an integrated product of paravectors, anyone can see a deep similarity between paravectors and vectors in Euclidean space. Using the integrated product we introduce geometric concepts into the algebra of paravectors, which contributes to its intuitive understanding.
\subsection{Parallelism and perpendicularity relations.}
\begin{definition}
Non-singular paravectors $\Gamma _{1}$ and $\Gamma _{2}$
are \textbf{parallel} $\left( \Gamma _{1}\parallel \Gamma
_{2}\right) $ if the vector product $( \Gamma_{1}, \Gamma_{2} \}=0$.
\end{definition}
\begin{theorem}
Two non-singular paravectors $\Gamma _{1}$ and $\Gamma_{2}$ are parallel if and only if there exists a number $\lambda \neq 0$ that \qquad $\Gamma _{1}=\lambda \Gamma _{2}$.
\end{theorem}
\begin{myproof}
Parallelism means that $\Gamma _{1}\Gamma_{2}^{-}=\alpha $, where $\alpha$ is a complex number. We multiply this equation on the right by paravector $\Gamma _{2}$

\vspace{1pt}$\Gamma _{1}(\Gamma _{2}^{-}\Gamma _{2})=\alpha \Gamma _{2}$

\noindent Since $\Gamma _{2}$ is non-singular then $\Gamma _{2}^{-}\Gamma_{2}=\beta$ is a non-zero complex number. Hence we get $\Gamma _{1}=\lambda \Gamma _{2}$, where $\lambda=\alpha /\beta $.
\end{myproof}
The above theorem shows that 
\begin{conclusion}
The parallelism satisfies the conditions of equivalence relation.
\end{conclusion}
\begin{definition}
Two non-singular paravectors $\Gamma _{1}$ and $\Gamma
_{2}$ are \textbf{perpendicular} $\left( \Gamma _{1}\perp \Gamma
_{2}\right) $ if the scalar product $\left\langle \Gamma_{1}, \Gamma_{2} \right\rangle=0$.
\end{definition}
The perpendicularity of paravectors has the same properties as the perpendicularity of vectors in Euclidean space.
\begin{theorem}
For each non-singular paravector:
\begin{enumerate}
\item \vspace{1pt}$\backsim \left( \Gamma \perp \Gamma \right) $
\item If $\Gamma _{1}\perp \Gamma _{2}$ then $\Gamma _{2}\perp \Gamma_{1}$
\item If $\Gamma _{1}\perp \Gamma _{2}$ and $\Gamma _{2}\parallel
\Gamma _{3}$ then $\Gamma _{1}\perp \Gamma _{3}$
\end{enumerate}
\end{theorem}
\begin{myproof}
\begin{enumerate}
\item \vspace{1pt}It follows by definition of perpendicularity.
\item It follows by the fact that the scalar product is symmetrical: 
$\left\langle \Gamma _{1},\Gamma _{2}\right\rangle =\left\langle \Gamma _{2},\Gamma _{1}\right\rangle$ and $\left\langle \Gamma _{1},\Gamma _{2}\right\rangle =0$
\item $\Gamma _{1}\perp \Gamma _{2}\iff \left( \Gamma _{1},\Gamma
_{2}\right\rangle =%
\begin{bmatrix}
0 \\ 
\pmb{\omega }%
\end{bmatrix}%
$
\newline
\vspace{1pt}\qquad \qquad $\Gamma _{2}\parallel \Gamma _{3}\iff \left( \Gamma
_{2},\Gamma _{3}\right\rangle =\lambda $
\vspace{1pt} hence \quad $\left( \Gamma _{1},\Gamma _{3}\right\rangle =\Gamma
_{1}\Gamma _{3}^{-}=\Gamma _{1}\Gamma _{2}^{-1}\Gamma _{2}\Gamma _{3}^{-}=%
\frac{\lambda }{\det \Gamma _{2}}\Gamma _{1}\Gamma _{2}^{-}=\begin{bmatrix}
0 \\ 
\frac{\lambda \pmb{\omega }}{\det \Gamma _{2}}\end{bmatrix}$
\end{enumerate}
\end{myproof}
\begin{conclusion}.
\begin{enumerate}
\item \vspace{1pt}
The paravector is perpendicular to itself if and only if it is singular.
\item \vspace{1pt}Paravectors mutually reversed (inversed) not to be parallel.
\item \vspace{1pt}
Orthogonal paravectors are parallel if and only if they are equal or opposite.
\end{enumerate}
\end{conclusion}
\begin{myproof}
\begin{enumerate}
\item It follows by assumption and definition of the singular paravector.
\item $\left( \Gamma,\Gamma^{-}\right\rangle =\Gamma \Gamma
=\begin{bmatrix}
\alpha \\ 
\pmb{\beta }\end{bmatrix}\begin{bmatrix}
\alpha \\ 
\pmb{\beta }
\end{bmatrix}=
\begin{bmatrix}
\alpha^{2}+\beta^{2} \\ 
2\alpha\pmb{\beta }\end{bmatrix}$
\item
Parallelism means that $\Lambda_{1} \Lambda_{2}^{-}=\lambda$, hence det$\Lambda_{1}$det$\Lambda_{2}=\lambda^2$. Paravectors are orthogonal hence $\lambda = \pm1$. First equality gives that $\Lambda_{1}=\Lambda_{2}$ or  $\Lambda_{1}=-\Lambda_{2}$.
\end{enumerate}
\end{myproof}
\begin{theorem}
For any non-singular paravectors $\Gamma _{1}$ and $\Gamma _{2}$ it occurs that:
\begin{enumerate}
\item If $\Gamma _{1}\perp \Gamma _{2}$ then $\Gamma _{1}^{\ast
}\perp \Gamma _{2}^{\ast }$

\item If $\Gamma _{1}\parallel \Gamma _{2}$ then $\Gamma _{1}^{\ast
}\parallel \Gamma _{2}^{\ast }$

\item If $\Gamma _{1}\parallel \Gamma _{2}$ then vig$\Gamma
_{1}\parallel$ vig$\Gamma _{2}$
\end{enumerate}
\end{theorem}
\begin{myproof}
\begin{enumerate}
\item $\left\langle \Gamma _{1},\Gamma _{2}\right\rangle =0\implies
\left\langle \Gamma _{2},\Gamma _{1}\right\rangle =0$, so $\left\langle \Gamma _{1}^{\ast },\Gamma _{2}^{\ast }\right\rangle
=\left\langle \Gamma _{2},\Gamma _{1}\right\rangle ^{\ast }=0$
\item $(\begin{bmatrix}
\alpha _{1} \\ 
\pmb{\beta }_{1}\end{bmatrix}\begin{bmatrix}
\alpha _{2} \\ 
-\pmb{\beta }_{2}
\end{bmatrix})_{V} =-\alpha _{1}\pmb{\beta }_{2}+\alpha _{2}\pmb{\beta }_{1}-i\pmb{\beta }_{1}\times \pmb{\beta }_{2}=\mathbf{0}$

Since complex vectors are governed by the same laws as real ones, so it must be:

$\alpha _{2}\pmb{\beta }_{1}-\alpha _{1}\pmb{\beta }_{2}=\pmb{
0\qquad }$and\qquad $i\pmb{\beta }_{1}\times \pmb{\beta }_{2}=\mathbf{0}$
\newline
On the other hand, we have $\left( \Gamma _{1}^{\ast },\Gamma _{2}^{\ast}\right\rangle =\left\langle \Gamma _{2},\Gamma _{1}\right) ^{\ast }=(\begin{bmatrix}
\alpha _{2} \\ 
-\pmb{\beta }_{2}
\end{bmatrix}
\begin{bmatrix}
\alpha _{1} \\ 
\pmb{\beta }_{1}
\end{bmatrix}
)^{\ast }$,
\newline
hence under the assumption

$(\begin{bmatrix}
\alpha _{2} \\ 
-\pmb{\beta }_{2}
\end{bmatrix}
\begin{bmatrix}
\alpha _{1} \\ 
\pmb{\beta }_{1}
\end{bmatrix})_{V}=
\alpha _{2}\pmb{\beta }_{1}-\alpha _{1}\pmb{\beta }_{2}-i\pmb{
\beta }_{2}\times \pmb{\beta }_{1}=\mathbf{0}$

\item $\left( \Gamma _{1}\Gamma _{1}^{\ast },\Gamma _{2}\Gamma _{2}^{\ast
}\right\rangle =\Gamma _{1}\Gamma _{1}^{\ast }(\Gamma _{2}\Gamma _{2}^{\ast
})^{-}=\Gamma _{1}\Gamma _{1}^{\ast }\Gamma _{2}^{\ast -}\Gamma
_{2}^{-}=\Gamma _{1}(\Gamma _{2}^{-}\Gamma _{1})^{\ast }\Gamma _{2}^{-}$

Under the assumption, the product in parentheses is a number, so we can move it in front of the product

$\lambda ^{\ast }\left( \Gamma _{1},\Gamma _{2}\right\rangle =\lambda ^{\ast
}\lambda $
\end{enumerate}
\end{myproof}

\begin{theorem}\label{th:3.4}
For any paravectors $\Gamma _{1}$ and $\Gamma _{2}$ polarization identity occurs:
\begin{equation}
\det (\Gamma _{1}+\Gamma _{2})=\det \Gamma _{1}+2\left\langle \Gamma
_{1},\Gamma _{2}\right\rangle +\det \Gamma _{2}
\end{equation}
\end{theorem}

\begin{myproof}
$\det (\Gamma _{1}+\Gamma _{2})=(\Gamma _{1}+\Gamma _{2})(\Gamma _{1}+\Gamma
_{2})^{-}=\Gamma _{1}\Gamma _{1}^{-}+\Gamma _{1}\Gamma _{2}^{-}+\Gamma
_{2}\Gamma _{1}^{-}+\Gamma _{2}\Gamma _{2}^{-}=$

$=\det \Gamma _{1}+\left( \Gamma _{1},\Gamma _{2}\right\rangle +\left(
\Gamma _{2},\Gamma _{1}\right\rangle +\det \Gamma _{2}=\det \Gamma
_{1}+2\left\langle \Gamma _{1},\Gamma _{2}\right\rangle +\det \Gamma _{2}$
\end{myproof}

\begin{conclusion}
Let paravectors $\Gamma_{1}$ and $\Gamma_{2}$ be perpendicular, then the determinant of these paravectors equals the sum of their determinants:
\begin{equation}
\det (\Gamma _{1}+\Gamma _{2})=\det \Gamma _{1}+\det \Gamma_{2} 
\end{equation}
\end{conclusion}

\vspace{1pt}This conclusion complies with the Pythagorean theorem in
Euclidean geometry.

Theorem 3.4 shows that paravectors meet the parallelogram law, which gives the structure constructed by us some of the characteristics of metric space.

\begin{conclusion}
For any paravectors $\Gamma _{1}$ and $\Gamma _{2}$ the parallelogram law occurs:
$$\text{det}\left( \Gamma_{1}+ \Gamma_{2} \right) +
\text{det}\left( \Gamma_{1}- \Gamma_{2} \right) =
2\text{det} \Gamma_{1} + 2\text{det} \Gamma_{2}$$
\end{conclusion}
\begin{definition}
Two paravectors
$\begin{bmatrix}
\alpha _{1} \\ 
\pmb{\beta }_{1}
\end{bmatrix}$ and $\begin{bmatrix}
\alpha _{2} \\ 
\pmb{\beta }_{2}
\end{bmatrix}$ are \textbf{spatially parallel} if $\pmb{\beta }%
_{1}\times \pmb{\beta }_{2}=\mathbf{0}$
\end{definition}

\begin{theorem}
If two non-singular paravectors are parallel, they are also spatially parallel.
\end{theorem}

\begin{myproof}
$\Gamma _{1}\parallel \Gamma _{2}\iff \alpha _{1}\pmb{\beta }_{2}-\alpha
_{2}\pmb{\beta }_{1}-i\pmb{\beta }_{1}\times \pmb{\beta }_{2}=0$

First we multiply the above equation by $\alpha _{1}\pmb{\beta }
_{2}$, and then by $\alpha _{2}\pmb{\beta }_{1}$. Hence we get two equations:

$\left( \alpha _{1}\pmb{\beta }_{2}\right) ^{2}-\alpha _{1}\alpha _{2}\pmb{\beta }_{1}\pmb{\beta }_{2}=0$

$\alpha _{1}\alpha _{2}\pmb{\beta }_{1}\pmb{\beta }_{2}-\left( \alpha
_{2}\pmb{\beta }_{1}\right) ^{2}=0$

The difference of the above equations yields $\left( \alpha _{1}\pmb{\beta }_{2}-\alpha _{2}
\pmb{\beta }_{1}\right) ^{2}=0$, which is true when $\alpha _{1}
\pmb{\beta }_{2}=\alpha _{2}\pmb{\beta }_{1}$.

Hence it follows that $\pmb{\beta }_{1}\times \pmb{\beta }_{2}=0$
\end{myproof}
As a consequence, we can say that

\begin{conclusion}
Spatial parallelism is an equivalence relation.
\end{conclusion}

\begin{conclusion}
The spatial parallelism of paravectors  is weaker than parallelism, ie.: If two paravectors are parallel, they must be spatially parallel, too, but in the opposite direction, implication no longer must occur.
\end{conclusion}

\begin{definition}
We call two singular paravectors $\Gamma _{1}$ and $\Gamma _{2}$  \textbf{singularly parallel} if $\left( \Gamma_{1},\Gamma _{2}\right\rangle =0$
\end{definition}

\begin{conclusion}
The following conclusions are easy to prove:
\begin{enumerate}
\item If two paravectors are singularly parallel, then they are singular.
\item Singular parallelism is an equivalence relation.
\end{enumerate}
\end{conclusion}
\subsection{Angles}
\begin{definition}
The \textbf{right angle} between two proper paravectors, denoted by $\angle ( \Gamma _{1},\Gamma _{2}\rangle $,  we call the paravector:
\begin{equation}
\Phi =\angle ( \Gamma _{1},\Gamma _{2}\rangle = \begin{bmatrix}
\text{cosi}\Phi \\ 
\text{\textbf{dex}}\Phi
\end{bmatrix}  \eqqcolon \frac{\left(
\Gamma _{1},\Gamma _{2}\right\rangle }{\left\vert \Gamma _{1}\right\vert
\left\vert \Gamma _{2}\right\vert },
\end{equation}
where we call the scalar component \textbf{cosinis},
and the vector component -- \textbf{dextis} of angle $\Phi$.
\end{definition}

\begin{definition}
The \textbf{left angle} between two proper paravectors we call the paravector:
\begin{equation}
\Phi =\angle \langle \Gamma _{1},\Gamma _{2}) = \begin{bmatrix}
\text{cosi}\Phi \\ 
\text{\textbf{sini}}\Phi
\end{bmatrix}  \eqqcolon \frac{\left\langle
\Gamma _{1},\Gamma _{2}\right) }{\left\vert \Gamma _{1}\right\vert
\left\vert \Gamma _{2}\right\vert },
\end{equation}
where we call the scalar component \textbf{cosinis},
and the vector component -- \textbf{sinis} of angle $\Phi$.
\end{definition}

\textbf{Note} - Please note that the components are not trigonometric (hyperbolic) functions - these are just names, given because angle components have the same properties as well-known trigonometric functions, making them easier to imagine and to remember.

These names are derived from Latin. \textit{Sinistram} means left, and \textit{dextram} means right.

In paravector space (which is not defined yet!) as well as in Euclidean space, we should identify a positive orientation. We don't do this in this paper, not to impose any restrictions. It seems that it will be necessary sooner or later, but not now.

\begin{definition}
We call an angle $\Phi =\Phi _{1}\Phi _{2}$ the \textbf{composition of} (left) \textbf{angles} $\Phi _{1}$ and 
$\Phi _{2}$.
\begin{align}
\Phi _{1}\Phi _{2}&=
\begin{bmatrix}
\text{cosi}\Phi _{1} \\ \nonumber
\text{\textbf{sini}}\Phi _{1}
\end{bmatrix}
\begin{bmatrix}
\text{cosi}\Phi _{2} \\ 
\text{\textbf{sini}}\Phi _{2}
\end{bmatrix}=\\
&=\begin{bmatrix}
\text{cosi}\Phi _{1}\text{cosi}\Phi _{2}+\text{\textbf{sini}}\Phi _{1}\text{\textbf{sini}}\Phi _{2} \\ \nonumber
\text{cosi}\Phi _{1}\text{\textbf{sini}}\Phi _{2}+\text{cosi}\Phi _{2}\text{\textbf{sini}}\Phi _{1}+i\text{\textbf{sini}}\Phi _{1}\times \text{\textbf{sini}}\Phi _{2}
\end{bmatrix}=\\
&=\begin{bmatrix}
\text{cosi}(\Phi _{1}\Phi _{2}) \\ \nonumber
\text{\textbf{sini}}(\Phi _{1}\Phi _{2})
\end{bmatrix}=
\begin{bmatrix}
\text{cosi}\Phi \\ \nonumber
\text{\textbf{sini}}\Phi
\end{bmatrix}=\Phi  
\end{align}
\end{definition}

We can write an analogous composition for the right angles. 

As a consequence, we can see further analogy with Euclidean trigonometry:

\begin{conclusion}The explement of an angle is:
\begin{eqnarray}
\nonumber \angle \langle \Gamma _{1},\Gamma _{2})^{-}
=\angle \langle \Gamma _{2},\Gamma _{1})& \qquad \text{for the left angle} \\
\text{or } \quad \angle ( \Gamma _{1},\Gamma _{2}\rangle^{-}
=\angle ( \Gamma _{2},\Gamma _{1}\rangle& \qquad \text{for the right angle},
\end{eqnarray}
which gives
\begin{itemize}
\item \vspace{1pt} cosi$\Phi ^{-}=$ cosi$\Phi $
\item \vspace{1pt}\textbf{sini}$\Phi ^{-}=$ $-$\textbf{sini}$\Phi $
\item \vspace{1pt}\textbf{dex}$\Phi ^{-}=$ $-$\textbf{dex}$\Phi $
\end{itemize}
\end{conclusion}

The left and right angles have an opposite orientation in space - they are not explementary.

$\angle\langle\Phi, \Gamma) = \frac{\Phi^{-} \Gamma}{\vert \Phi \vert \vert \Gamma \vert}$, 
$\qquad \angle (\Phi, \Gamma\rangle = \frac{\Phi \Gamma^{-}}{\vert \Phi \vert \vert \Gamma \vert}, \qquad$ hence the composition of this angles is:
$$\angle\langle \Phi, \Gamma)\angle (\Phi, \Gamma\rangle=
\frac{\Phi^{-} \Gamma \Phi \Gamma^{-}}{\text{det}(\Phi \Gamma)}\neq 1$$
The explement of the left angle $\angle\langle\Phi, \Gamma)$ is the left angle $\angle\langle \Gamma,\Phi)$, and the same occurs for the right angle.

In Table 3.1 are shown the properties of left angle components to simply and intuitively justify the names chosen for them. The right angle has analogous formulas.

\vspace{5pt}Tab. 3.1\qquad General recurrence formulas for components of the left angle.

\begin{tabular}{|l|l|}
\hline
$\Phi $ - Orthogonal paravector & $\text{cosi}\Phi \in C$ and $\text{\textbf{sini}}\Phi \in C^{3}$ \\ \hline\hline
Determinant of $\Phi $ & cosi$^{2}\Phi $ - $\text{sini}^{2}\Phi =1$ \\ \hline
Doubling & $\text{cosi}(\Phi ^{2})=$ cosi$^{2}\Phi $ + $\text{sini}%
^{2}\Phi $ \\ 
the angle & $\mathbf{sini}(\Phi ^{2})=$ 2cosi$\Phi\mathbf{sini}\Phi $ \\ 
\hline
Angles & $\text{cosi}(\Phi _{1}\Phi _{2})=\text{cosi}\Phi _{1}\text{%
cosi}\Phi _{2}+\text{\textbf{sini}}\Phi _{1}\text{\textbf{sini}}\Phi _{2}$
\\ 
composition & $\text{\textbf{sini}}(\Phi _{1}\Phi _{2})=\text{cosi}\Phi
_{1}\text{\textbf{sini}}\Phi _{2}+\text{cosi}\Phi _{2}\text{\textbf{sini}}
\Phi _{1}+i\text{\textbf{sini}}\Phi _{1}\times \text{\textbf{sini}}\Phi _{2}$
\\ \hline
\multirow{2}{*}{Angle explement}& cosi$\Phi ^{-}=$ cosi$\Phi $ \\ 
 & \textbf{sini}$\Phi ^{-}=$ $-$\textbf{sini}$\Phi $ \\ \hline
\end{tabular}
\vspace{5pt}
\begin{conclusion}
If $Im($cosi$\Phi )=0$ and $Im($\textbf{sini}$\Phi )=\mathbf{0}$, then the nature of the angle $\Phi $ is hyperbolic. If $Im($cosi$\Phi )=0$ and $Re($\textbf{sini}$\Phi )=\mathbf{0}$, then the nature of the angle is trigonometric, which is shown by the determinant of this angle.
\end{conclusion}

\subsection{Similarity and rotation}

\begin{definition} \label{def:2.3.7}
Two paravectors $\Gamma_1$ and $\Gamma_2$ are \textbf{similar} if there exists a non-singular paravector $\Phi$, that $\Phi \Gamma_{1} = \Gamma_{2} \Phi$. We call the paravector $\Phi$ \textbf{axis of similarity}.
\end{definition}
Similarity can be shown in any other way $\Gamma^{\prime} = \Phi^{- 1} \Gamma \Phi$, which is like many authors define the rotation. We would like to associate the rotation with the angle between the rotated paravector and its image after the turning unequivocally. If the paravector $\Phi$ is improper, then it's impossible to determine this angle. For this reason, it was necessary to clarify the definition of rotation.           
          	\begin{definition}
A similarity whose axis is a proper paravector is called the \textbf{rotation}.
            \end{definition}             
We need to orient the rotation in accordance with the angles, therefore we will distinguish left and right rotations. The left rotation will take the form $\Gamma^{\prime }=\Phi^{-1}\Gamma \Phi$, and the right: $\Gamma^{\prime }=\Phi\Gamma \Phi^{-1}$. Since the paravector $\Phi$ is proper, then the module $\vert \Phi \vert$ exists. So, we can exhibit the rotation in the form:
\begin{equation}
\Gamma ^{\prime }=\Lambda^{-}\Gamma \Lambda
\end{equation}
       
\noindent  where paravector $\Lambda = \frac{\Phi}{\vert \Phi \vert} $ is the axis of rotation and also determines the value of rotation. The axis of rotation in space is determined by the spatial component of paravector $\Lambda$.
                            
The properties of rotations on the space built by us are so general that we cannot restrict them to the rotation in Euclidean sense only. In cases where cosi$\Lambda $ is a real number, and \textbf{sini}$\Lambda $ is an imaginary vector (or an angle $\Lambda $ is a special paravector), we deal with an Euclidean rotation. For the paravector $\Lambda = \begin{bmatrix}
                    \cos \varphi \\ 
                    i\textbf{n}\sin \varphi
                \end{bmatrix} $ we have spatial rotation by the angle $2\varphi$ about the axis defined by vector \textbf{n}. Here $\varphi$ is a traditional angle, and sine/cosine are trigonometric functions. The fact that the paravector angle (despite the similarities) is something other than the Euclidean angle we can see by examining the right angle between $\Gamma$ and its rotated image $\Gamma^{\prime}=\Lambda^{-}\Gamma \Lambda$. The right angle between paravectors $\Gamma$ and $\Gamma^{\prime}$ can be written
                
\begin{equation}
 \angle (\Gamma,\Gamma^{\prime}\rangle = 
                \frac{\Gamma(\Lambda^{-}\Gamma \Lambda)^{-}}{\text{det}\Gamma}=
                \frac{(\Gamma\Lambda^{-})(\Gamma^{-} \Lambda)}{\text{det}\Gamma}=
                \angle(\Gamma,\Lambda\rangle \angle\langle\Gamma,\Lambda)
\end{equation}

\begin{conclusion}                				
In the space the angle between a paravector and its image after turning, is a combination of the right and left angles between these paravectors and the axis of rotation. 
\end{conclusion}
In the case when the axis of rotation is a real paravector, W. Baylis says that such a rotation is a Lorentz transformation of the electric field.

Below we show the obvious properties of similarity:
\begin{theorem}\label{th:3.6} 
For each similar paravector:
\begin{enumerate}
\item Similar paravectors must have the same scalar components. In other words, similarity is a spatial relationship.
\item If any paravectors are spatially parallel to the axis of similarity and they are similar, then they are identical.
\item Parallel axes represent the same similarity.
\item Similarity is an equivalence relation.
\end{enumerate}
\end{theorem}

\begin{myproof}
\begin{enumerate}
\item $\Phi ^{-1}\Gamma \Phi =\frac{1}{\alpha ^{2}-\beta ^{2}}
\begin{bmatrix}
\alpha \\ 
\pmb{\beta }
\end{bmatrix}
\begin{bmatrix}
\tau \\ 
\pmb{\varpi }
\end{bmatrix}
\begin{bmatrix}
\alpha \\ 
-\pmb{\beta }
\end{bmatrix}
=\frac{1}{\alpha ^{2}-\beta ^{2}}
\begin{bmatrix}
\left( \alpha ^{2}-\beta ^{2}\right) \tau \\ 
\left( \alpha ^{2}-\beta ^{2}\right) \pmb{\varpi -}2\left[ i\pmb{\beta 
}\times \pmb{\varpi }+\pmb{\beta }\times \left( \pmb{\beta }\times 
\pmb{\varpi }\right) \right]
\end{bmatrix}$
\item It follows from the previous property and $\pmb{\beta }\times \pmb{\varpi} =\mathbf{0}$
\item Let $\Phi _{1}=\lambda \Phi _{2}$

$\frac{1}{\det \Phi _{1}}\Phi _{1}^{-}\Gamma \Phi _{1}=\frac{1}{\lambda
^{2}\det \Phi _{2}}\lambda ^{2}\Phi _{2}^{-}\Gamma \Phi _{2}=\frac{1}{\det
\Phi _{2}}\Phi _{2}^{-}\Gamma \Phi _{2}$
\item The proof is simple, so we leave it to the reader.
\end{enumerate}
\end{myproof}
By the theorem \ref{th:3.6} we draw the following conclusions for rotations:

\begin{conclusion} \label{con:3.10}
For each rotated paravector:
\begin{enumerate}
\item Rotation does not change the scalar component. In other words, rotation is a spatial relationship.                  
\item Rotation does not change the paravector which is spatially parallel to the axis of this rotation.				
\item Parallel axes represent the same rotation.
\end{enumerate}
\end{conclusion}

\noindent Using paravectors, we can easily introduce Euler angles, or the composition of angles on the planes with normal 
$\mathbf{n}_{1}$ and $\mathbf{n}_{2}$
\begin{equation}
\begin{bmatrix}
\cos \varphi \\ 
i\mathbf{n}\sin \varphi
\end{bmatrix}
=
\begin{bmatrix}
\cos \varphi _{1} \\ 
i\mathbf{n}_{1}\sin \varphi_{1}
\end{bmatrix}
\begin{bmatrix}
\cos \varphi_{2} \\ 
i\mathbf{n}_{2}\sin \varphi_{2}
\end{bmatrix}
\qquad \qquad \left( \left\vert \mathbf{n}_{i}\right\vert =1\right)
\end{equation}

\vspace{1pt}Please note that the above angles $\varphi$ and functions have a trigonometric sense in real Euclidean space!

Any vector in the space can be written as a paravector $\begin{bmatrix}
									0 \\ 
									i\textbf{w}
								\end{bmatrix}$ (imaginary vector, because if it were real, then the paravector would not have a module). The angle between two vectors $\textbf{w}_1$ and $\textbf{w}_2$, is
                                
\begin{equation}
\frac{1}{\vert \textbf{w}_{1} \vert} \begin{bmatrix}
									0 \\ 
									i\textbf{w}_{1}
								\end{bmatrix}
                                \frac{1}{\vert \textbf{w}_{2} \vert} \begin{bmatrix}
									0 \\ 
									-i\textbf{w}_{2}
								\end{bmatrix}=
                                \begin{bmatrix}
									\frac{\textbf{w}_{1}\textbf{w}_{2}}
                                    {\vert \textbf{w}_{1} \vert\vert \textbf{w}_{2} \vert} \\ 
									i\frac{\textbf{w}_{1} \times \textbf{w}_{2}}
                                    {\vert \textbf{w}_{1} \vert\vert \textbf{w}_{2} \vert}
								\end{bmatrix}=
                                \begin{bmatrix}
                    \cos x \\ 
                    i\textbf{n}\sin x
                \end{bmatrix}
\end{equation}

\begin{definition}
\textbf{Mirror symmetry} with respect to the plane of normal $\mathbf{n}$ is called conversion such that:
\begin{equation}
\begin{bmatrix}
0 \\ 
i\mathbf{w}^{\prime}
\end{bmatrix}=
\begin{bmatrix}
0 \\ 
i\mathbf{n}
\end{bmatrix}
\begin{bmatrix}
0 \\ 
i\mathbf{w}
\end{bmatrix}
\begin{bmatrix}
0 \\ 
i\mathbf{n}
\end{bmatrix}
\end{equation}
\end{definition}

\vspace{1pt}Let $\mathbf{n}$ be a vector normal to the plane N. Take any non-zero vector $\mathbf{w}$. This vector can be decomposed into a vector parallel to $\mathbf{n}$ and perpendicular (orthogonal projection of vector $\mathbf{w}$ on the plane N), respectively:

$\mathbf{n}\left( \mathbf{wn}\right) \qquad $and\qquad $\mathbf{w}-\mathbf{n}
\left( \mathbf{wn}\right) =\mathbf{n}^{2}\mathbf{w}-\mathbf{n}\left( \mathbf{wn}\right) =\left( \mathbf{n}\times \mathbf{w}\right) \times \mathbf{n},$

hence vector $\mathbf{w}=\mathbf{n}\left( \mathbf{wn}\right) +\left( 
\mathbf{n}\times \mathbf{w}\right) \times \mathbf{n}.$

\vspace{1pt}Mirror symmetry changes the sign of the component perpendicular to the plane of symmetry (parallel to \textbf{n}), but when we take a non-zero scalar, we get:

\[
\begin{bmatrix}
0 \\ 
i\mathbf{n}
\end{bmatrix}
\begin{bmatrix}
a \\ 
\mathbf{w}
\end{bmatrix}
\begin{bmatrix}
0 \\ 
i\mathbf{n}
\end{bmatrix}=
\begin{bmatrix}
-a \\ 
-\mathbf{n}\left( \mathbf{wn}\right) +\left( \mathbf{n}\times \mathbf{w}
\right) \times \mathbf{n}
\end{bmatrix}
\]

Hence we see that mirror symmetry is not similarity in the meaning of definition \ref{def:2.3.7}, because it changes the sign of a scalar.

\vspace{1pt}Mirror symmetry can be generalized to complex paravectors:

\begin{equation}
\frac{1}{-\omega ^{2}}%
\begin{bmatrix}
0 \\ 
\pmb{\omega }
\end{bmatrix}
\begin{bmatrix}
\alpha \\ 
\pmb{\beta }
\end{bmatrix}
\begin{bmatrix}
0 \\ 
\pmb{\omega }
\end{bmatrix}
=\begin{bmatrix}
-\alpha \\ 
\omega ^{-2}[-\pmb{\omega }\left( \pmb{\beta \omega }\right) +\left( 
\pmb{\omega }\times \pmb{\beta }\right) \times \pmb{\omega ]}
\end{bmatrix}
\end{equation}

\vspace{1pt}As was to be expected, rotation can be presented in the form of a composition of two mirror symmetries. The paravector parallel to the both planes of symmetry sets the axis of rotation.
\[
\frac{1}{\omega _{1}^{2}\omega _{2}^{2}}
\begin{bmatrix}
0 \\ 
\pmb{\omega }_{2}
\end{bmatrix}
\begin{bmatrix}
0 \\ 
\pmb{\omega }_{1}
\end{bmatrix}
\begin{bmatrix}
\alpha \\ 
\pmb{\beta }
\end{bmatrix}
\begin{bmatrix}
0 \\ 
\pmb{\omega }_{1}
\end{bmatrix}
\begin{bmatrix}
0 \\ 
\pmb{\omega }_{2}
\end{bmatrix}
=\frac{1}{\left( \pmb{\omega }_{1}\pmb{\omega }
_{2}\right) ^{2}+\left( \pmb{\omega }_{1}\times \pmb{\omega }
_{2}\right) ^{2}}
\begin{bmatrix}
\pmb{\omega }_{1}\pmb{\omega }_{2} \\ 
-i\pmb{\omega }_{1}\times \pmb{\omega }_{2}
\end{bmatrix}
\begin{bmatrix}
\alpha \\ 
\pmb{\beta }
\end{bmatrix}
\begin{bmatrix}
\pmb{\omega }_{1}\pmb{\omega }_{2} \\ 
i\pmb{\omega }_{1}\times \pmb{\omega }_{2}
\end{bmatrix}
\]

\vspace{1pt}\textbf{Axial symmetry} is nothing else but a straight angle rotation around the vector $\pmb{\omega }$

\begin{equation}
\frac{1}{\omega ^{2}}
\begin{bmatrix}
0 \\ 
-i\pmb{\omega }
\end{bmatrix}
\begin{bmatrix}
\alpha \\ 
\pmb{\beta }
\end{bmatrix}
\begin{bmatrix}
0 \\ 
i\pmb{\omega }
\end{bmatrix}
=\begin{bmatrix}
\alpha \\ 
\omega ^{-2}[\pmb{\omega }\left( \pmb{\beta \omega }\right) -\left( 
\pmb{\omega }\times \pmb{\beta }\right) \times \pmb{\omega ]}
\end{bmatrix}
\end{equation}

From the above discussion we can see that paravectors, even though are of complex construction and have no vector metric, have geometrical features of vectors so that they become imaginable.

\section{Matrix representation of paravectors}

Based on the definition of paravectors multiplication \eqref{def:1.6}, the equation
\begin{equation}
X_{2}=\Gamma X_{1} =
\begin{bmatrix}
\alpha_{2} \\ 
\pmb{\beta }_{2}
\end{bmatrix} =
\begin{bmatrix}
\alpha \\ 
\pmb{\beta }
\end{bmatrix}
\begin{bmatrix}
\alpha_{1} \\ 
\pmb{\beta }_{1}
\end{bmatrix}
\end{equation}
we can denote
\begin{equation}
\begin{bmatrix}
\alpha_{2} \\ 
\pmb{\beta }_{2}
\end{bmatrix} =
\begin{bmatrix}
\alpha \alpha _{1}+\pmb{\beta }\pmb{\beta }_{1} \\ 
\alpha\pmb{\beta }_{1}+\alpha _{1}\pmb{\beta } +i\pmb{\beta 
}\times \pmb{\beta }_{1}
\end{bmatrix}
\end{equation}
The above equation is a system of linear equations, which can be exhibit in matrix form
\begin{equation}
\begin{bmatrix}
\alpha_{2} \\ 
\beta_{2x}\\
\beta_{2y}\\
\beta_{2z}
\end{bmatrix} =
\begin{bmatrix}
\alpha & \beta _{x} & \beta _{y} & \beta _{z} \\ 
\beta _{x} & \alpha & -i\beta _{z} & i\beta _{y} \\ 
\beta _{y} & i\beta _{z} & \alpha & -i\beta _{x} \\ 
\beta _{z} & -i\beta _{y} & i\beta _{x} & \alpha
\end{bmatrix}
\begin{bmatrix}
\alpha_{1} \\ 
\beta_{1x}\\
\beta_{1y}\\
\beta_{1z}
\end{bmatrix}
\end{equation}

Anyone can see that above equation is equivalent to
\begin{equation}
\begin{bmatrix}
\alpha_{2} & \beta _{2x} & \beta _{2y} & \beta _{2z} \\ 
\beta _{2x} & \alpha_{2} & -i\beta _{2z} & i\beta _{2y} \\ 
\beta _{2y} & i\beta _{2z} & \alpha_{2} & -i\beta _{2x} \\ 
\beta _{2z} & -i\beta _{2y} & i\beta _{2x} & \alpha_{2}
\end{bmatrix} =
\begin{bmatrix}
\alpha & \beta _{x} & \beta _{y} & \beta _{z} \\ 
\beta _{x} & \alpha & -i\beta _{z} & i\beta _{y} \\ 
\beta _{y} & i\beta _{z} & \alpha & -i\beta _{x} \\ 
\beta _{z} & -i\beta _{y} & i\beta _{x} & \alpha
\end{bmatrix}
\begin{bmatrix}
\alpha_{1} & \beta _{1x} & \beta _{1y} & \beta _{1z} \\ 
\beta _{1x} & \alpha_{1} & -i\beta _{1z} & i\beta _{1y} \\ 
\beta _{1y} & i\beta _{1z} & \alpha_{1} & -i\beta _{1x} \\ 
\beta _{1z} & -i\beta _{1y} & i\beta _{1x} & \alpha_{1}
\end{bmatrix}
\end{equation}

Therefore, each paravector $\Gamma =
\begin{bmatrix}
\alpha \\ 
\pmb{\beta }
\end{bmatrix}$  is equivalent to a matrix

\begin{equation}
\begin{bmatrix}
\alpha & \beta _{x} & \beta _{y} & \beta _{z} \\ 
\beta _{x} & \alpha & -i\beta _{z} & i\beta _{y} \\ 
\beta _{y} & i\beta _{z} & \alpha & -i\beta _{x} \\ 
\beta _{z} & -i\beta _{y} & i\beta _{x} & \alpha
\end{bmatrix}
\end{equation}

The determinant of the above matrix is $\label{eq:2.30} \left( \alpha ^{2}-\beta
^{2}\right) ^{2}=\left( \Gamma \Gamma ^{-}\right) ^{2}=\left( \det \Gamma
\right) ^{2},$ hence paravector $\Gamma ^{-1}=\frac{\Gamma ^{-}}{\det \Gamma }$ corresponds to the matrix inverse of the above one.

Since the inverse paravector should correspond to the transposed matrix, then we were considering naming it a transposed paravector. But this transposition is not complete because the first row and first column  are not subject to transposition:
\begin{equation}
\Gamma \longleftrightarrow
\begin{bmatrix}
\alpha & \beta _{x} & \beta _{y} & \beta _{z} \\ 
\beta _{x} & \alpha & -i\beta _{z} & i\beta _{y} \\ 
\beta _{y} & i\beta _{z} & \alpha & -i\beta _{x} \\ 
\beta _{z} & -i\beta _{y} & i\beta _{x} & \alpha
\end{bmatrix}%
\qquad \qquad \qquad \qquad \Gamma ^{-}\longleftrightarrow%
\begin{bmatrix}
\alpha & -\beta _{x} & -\beta _{y} & -\beta _{z} \\ 
-\beta _{x} & \alpha & i\beta _{z} & -i\beta _{y} \\ 
-\beta _{y} & -i\beta _{z} & \alpha & i\beta _{x} \\ 
-\beta _{z} & i\beta _{y} & -i\beta _{x} & \alpha%
\end{bmatrix}%
\end{equation}

The geometric meaning of this paravector corresponds to the reverse direction in space, so it was decided to leave the name: reverse paravector.

\begin{conclusion}Some of the matrices counterparts:
\begin{enumerate}
\item The singular paravector corresponds to the singular matrix.
\item The conjugate paravector corresponds to the Hermitian conjugate matrix.
\end{enumerate}
\end{conclusion}

\begin{proof} of the 2nd point.

\vspace{10pt}
$
\begin{bmatrix}
\alpha \\ 
\pmb{\beta }
\end{bmatrix}^{\ast}\longleftrightarrow
\begin{bmatrix}
\alpha ^{\ast } & \beta _{x}^{\ast } & \beta _{y}^{\ast } & \beta _{z}^{\ast
} \\ 
\beta _{x}^{\ast } & \alpha ^{\ast } & -i\beta _{z}^{\ast } & i\beta
_{y}^{\ast } \\ 
\beta _{y}^{\ast } & i\beta _{z}^{\ast } & \alpha ^{\ast } & -i\beta
_{x}^{\ast } \\ 
\beta _{z}^{\ast } & -i\beta _{y}^{\ast } & i\beta _{x}^{\ast } & \alpha
^{\ast }
\end{bmatrix}
=
\begin{bmatrix}
a-id & b_{x}-ic_{x} & b_{y}-ic_{y} & b_{z}-ic_{z} \\ 
b_{x}-ic_{x} & a-id & -ib_{z}-c_{z} & ib_{y}+c_{y} \\ 
b_{y}-ic_{y} & ib_{z}+c_{z} & a-id & -ib_{x}-c_{x} \\ 
b_{z}-ic_{z} & -ib_{y}-c_{y} & ib_{x}+c_{x} & a-id
\end{bmatrix}=$

$=\begin{bmatrix}
a+id & b_{x}+ic_{x} & b_{y}+ic_{y} & b_{z}+ic_{z} \\ 
b_{x}+ic_{x} & a+id & ib_{z}-c_{z} & -ib_{y}+c_{y} \\ 
b_{y}+ic_{y} & -ib_{z}+c_{z} & a+id & ib_{x}-c_{x} \\ 
b_{z}+ic_{z} & ib_{y}-c_{y} & -ib_{x}+c_{x} & a+id
\end{bmatrix}
^{\ast }=$

$=\begin{bmatrix}
a+id & b_{x}+ic_{x} & b_{y}+ic_{y} & b_{z}+ic_{z} \\ 
b_{x}+ic_{x} & a+id & i(b_{z}+ic_{z}) & -i(b_{y}+ic_{y}) \\ 
b_{y}+ic_{y} & -i(b_{z}+ic_{z}) & a+id & i(b_{x}+ic_{x}) \\ 
b_{z}+ic_{z} & i(b_{y}+ic_{y}) & -i(b_{x}+ic_{x}) & a+id%
\end{bmatrix}
^{\ast }=
\begin{bmatrix}
\alpha & \beta _{x} & \beta _{y} & \beta _{z} \\ 
\beta _{x} & \alpha & -i\beta _{z} & i\beta _{y} \\ 
\beta _{y} & i\beta _{z} & \alpha & -i\beta _{x} \\ 
\beta _{z} & -i\beta _{y} & i\beta _{x} & \alpha%
\end{bmatrix}
^{\ast T}$
\end{proof}

\section{Orthogonal transformations}

\begin{definition}
Linear transformation represented by a paravector will be called a \textbf{paravector transformation}.
\end{definition}

\begin{definition}
A paravector transformation is called \textbf{orthogonal} if it preserves the scalar product of paravectors.
\end{definition}
Or equivalently: A paravector transformation is called orthogonal if its determinant is equal to 1.

From the condition of orthogonal transformation it can be seen that its paravector must be proper (def.\ref{def:2.1.9}), and thus it has a module (def.\ref{def:2.1.12}).

\begin{conclusion}
An orthogonal transformation is represented by paravector:
\begin{equation}
\Lambda =\frac{1}{\sqrt{\alpha ^{2}-\beta ^{2}}}
\begin{bmatrix}
\alpha \\ 
\pmb{\beta }
\end{bmatrix}
=\frac{1}{\sqrt{a^{2}-b^{2}+c^{2}-d^{2}}}
\begin{bmatrix}
a+id \\ 
\mathbf{b+}i\mathbf{c}
\end{bmatrix}
\end{equation}
such that $ad=\mathbf{bc}$
\end{conclusion}

\begin{definition}
A transformation which preserves determinants is called an \textbf{isometric transformation}.
\end{definition}

\begin{conclusion}
An orthogonal transformation is isometric.
\end{conclusion}

\begin{example}
Invariance of the sphere under orthogonal transformation.

An equation of the sphere of radius $r$ can be written:
\begin{equation}
r^{2}-x^{2}=
\begin{bmatrix}
r \\ 
\mathbf{x}
\end{bmatrix}
\begin{bmatrix}
r \\ 
-\mathbf{x}
\end{bmatrix}
=
\begin{bmatrix}
r \\ 
\mathbf{x}
\end{bmatrix}
\begin{bmatrix}
r \\ 
\mathbf{x}
\end{bmatrix}
^{-}=0
\end{equation}
Since the sphere equation is the determinant of a singular paravector $\begin{bmatrix}
						r \\ 
						\mathbf{x}
					\end{bmatrix}$, so based on conclusion 5.2 we can see that in the complex space the spherical shape must be invariant with respect to the discussed transformation.
\end{example}

\begin{theorem}
If paravectors $\Gamma _{1}$ and $\Gamma _{2}$ are parallel, and $\Gamma ^{\prime }=\Lambda \Gamma $ or $\Gamma ^{\prime }=\Gamma\Lambda$ (where $\det \Lambda =1$),
then paravectors vig$\Gamma _{1}^{\prime }$ and vig$\Gamma _{2}^{\prime }$ are parallel, too.
\end{theorem}

So, a paravector transformation preserves the parallelism of vigors of these paravectors.

\begin{theorem}
Let $\Lambda $ be a orthogonal paravector, then
\begin{enumerate}
\item Transformation $\Gamma ^{\prime }=\Gamma \Lambda $ preserves a scalar product of the products $\Gamma ^{\ast }\Gamma $,
\item Transformation $\Gamma ^{\prime }=\Lambda \Gamma $ preserves a scalar product of the vigors $\Gamma \Gamma ^{\ast }$.
\end{enumerate}
\end{theorem}

\begin{myproof}

\begin{enumerate}
\item                            	
   Let $\Gamma^{\prime} = \Gamma \Lambda$ and det$\Lambda = 1$
                                
    $\Gamma^{\prime \ast}_{1} \Gamma^{\prime}_{1} 
                                (\Gamma^{\prime \ast}_{2} \Gamma^{\prime}_{2})^{-}
                                =(\Gamma_{1} \Lambda)^{\ast} \Gamma_{1} \Lambda 
                                [(\Gamma_{2} \Lambda)^{\ast} \Gamma_{2} \Lambda]^{-}=$
                                
							$=\Lambda^{\ast} \Gamma_{1}^{\ast}\Gamma_{1}\Lambda\Lambda^{-}
                            \Gamma_{2}^{-}\Gamma_{2}^{\ast -}\Lambda^{\ast -}=$
                            
                            $=\Lambda^{\ast}(\Gamma_{1}^{\ast}\Gamma_{1})
                            (\Gamma_{2}^{\ast}\Gamma_{2})^{-}\Lambda^{\ast -}$
                            
  which completes the proof, since from the conclusion \ref{con:3.10}.1, a rotation does not change the scalar.

\item Proof runs the same way as above. Let change the transformation $\Gamma^{\prime} = \Lambda\Gamma $
 
                                $\Gamma^{\prime}_{1} \Gamma^{\prime \ast}_{1}
                                ( \Gamma^{\prime}_{2}\Gamma^{\prime \ast}_{2})^{-}
                                =\Lambda\Gamma_{1}(\Lambda\Gamma_{1})^{\ast}  
                                [\Lambda\Gamma_{2}(\Lambda\Gamma_{2})^{\ast}]^{-}=$
                                
                            $=\Lambda\Gamma_{1}\Gamma_{1}^{\ast}\Lambda^{\ast}
                            \Lambda^{\ast -}\Gamma_{2}^{\ast -}\Gamma_{2}^{-}\Lambda^{-}=$
                            
                            $=\Lambda(\Gamma_{1}\Gamma_{1}^{\ast})
                        (\Gamma_{2}\Gamma_{2}^{\ast})^{-}\Lambda^{-}$
\end{enumerate}
\end{myproof}

\section{Summary}

In school courses physical sizes are always divided into scalars or vectors. Their natural generalization are paravectors which have both characteristics of integers, as well as the geometric properties of vectors. We have tried to present paravectors in order to show their similarity to vectors in Euclidean space.

Since paravectors together with the operation of summation form an Abelian group and there the scalar product of paravectors is defined we can say that paravectors form an unitary space over the complex numbers field. However, there is a significant difference between the definition of the scalar product of paravectors and definitions commonly known, becouse the product of any paravector with itself cannot be a real number. The conclusion is that space of paravectors is not normed, but we can define the function wich possesses, some of properties the same as the square of norm has. This function is a determinant which fulfils the parallelogram law and polarization equation. Unfortunately, the determinant is not enough to introduce the ordering relation, because their values of are complex numbers and complex numbers are not ordered. The concept of a norm can be only entered for proper and singular paravectors (a module of paravector). The trouble is that proper paravectors together with the addition operation do not form a group. As can be seen, the issue is so wide that its solution will be presented in some of the next publications.

By exploring various properties of paravectors we have found that different groups of them have different properties, in spite of the same construction. Some of them act as vectors, and other as matrices, therefore:
                
\begin{itemize}
\item	Additive paravectors (i.e. coordinates or field functions), denoted as bellow, are called  \textbf{four-vectors} traditionally.

\begin{equation}
\mathbb{X}  \eqqcolon
                            \begin{pmatrix}
                            \Delta t \\ 
                            \Delta \mathbf{x}%
                            \end{pmatrix}%
\end{equation}						
                    
\item  Paravectors which are not additive (transformation parameters i.e. speed or rotation) are denoted in square brackets:

\begin{equation}
               \Gamma  \eqqcolon
                \begin{bmatrix}
                \alpha \\ 
                \boldsymbol{\beta }
                \end{bmatrix}
                =
                \begin{bmatrix}
                a+id \\ 
                \textbf{b}+i\textbf{c}
                \end{bmatrix}
\end{equation}            	
            
 \end{itemize}
Interesting and surprising results of applications of the presented paravectors in field theory will be shown in a series of next articles. We hope that the notation proposed by us will replace in future the currently applied formalisms in algebra of physical space (APS) and will give a new perspective to look at the space-time.

\newpage
\section*{Appendix: Reference to the existing formalisms}

\vspace{1pt}In Clifford Algebra both the formalism and the terminology used by different authors are very diverse. The formalism based on Grassmann algebra, promoted by leading contemporary researchers is dominant now. We based on the matrix algebra due to:
\begin{itemize}
\item reputation of the matrix calculus. Grassmann algebra, although belongs to the  mathematical classics, is skillfully used by a relatively small and tight group of researchers, while the matrix algebra is a commonly used mathematical tool.
\item accounting transparency. Implementation of the activities in which three or more paravectors are  multiplying using the concepts of multivector is very complicated.
\end{itemize}

Paravectors and complex quaternions have the same construction but operations of multiplication are defined differently for them. In quaternion multiplication there is no imaginary one at the cross product. It gives that the multiplication of complex quaternions is not associative unlike the multiplication of paravectors. Forty years ago the notation similar to ours was used by Aharonow, Farach and Poole \cite{Aharonow} - \cite{Aharonow_Vector}, and recently W.E.Baylis \cite{Baylis_Algebra}.
We differ from them in the fact that the pair scalar-vector in the said articles is written in a row, and we have a column. The row notation is more practical for edition, while the column notation seems favorable for calculation, because it creates separation transparency into the scalar and vector portions and provides it a quick and easy way to transition to the matrix representation.

\vspace{1pt}For the reader familiar with the classic paravectors, it can be important to compare our formalism with a classic one. Since the closest to our notation seems the notation applied by prof. William Baylis and we often refer to his works, we give a concise relationship between our and \href{http://cronus.uwindsor.ca/baylis-research}{his} definitions.  (\href{http://en.wikipedia.org/wiki/Paravector}{Wiki})

\begin{tabular}[h]{|l|l||l|l|}
\hline
\multicolumn{2}{|c||}{W.E.Baylis} & \multicolumn{2}{c|}{This article} \\ \hline \hline
 paravector & $q = a + id + \mathbf{b} + i\mathbf{c}$  & paravector  & $\Gamma=	\begin{bmatrix}
					a+id \\ 
					\mathbf{b}+i\mathbf{c}
					\end{bmatrix}$         \\ \hline
 bar conjugation $^{\ast )}$ & $\overline{q}= a + id - \mathbf{b}- i\mathbf{c}$ & reversion  & $\Gamma ^{-}=
					\begin{bmatrix}
					a+id \\ 
					-\mathbf{b}-i\mathbf{c}
					\end{bmatrix}$                  \\ \hline
 Hermitean conjugation & $q^{\dagger }= a - id + \mathbf{b} - i\mathbf{c}$ & conjugation & $\Gamma ^{\ast }=\begin{bmatrix}
					a-id \\ 
					\textbf{b}-i\textbf{c}
					\end{bmatrix}$ \\ \hline
\makecell[l]{gradient operator \\ (or paragradient)} & $\partial =\partial /\partial t-\nabla $  & \makecell[l]{reversed \\ differential operator\\(or 4-gradient)}  & $\partial ^{-}=\begin{bmatrix}
					\partial /\partial t \\ 
					-\nabla 
					\end{bmatrix}$   \\ \hline
 \makecell[l]{spatially reversed \\ gradient operator} & $\partial =\partial /\partial t+\nabla $  & \makecell[l]{differential operator\\(or 4-divergence)}  & $\partial =\begin{bmatrix}
					\partial /\partial t \\ 
					\nabla 
					\end{bmatrix}$ \\ \hline
  &  &   determinant & $\det $ $\Gamma =\Gamma \Gamma ^{-}$                   \\ \hline
  &  &   module  $^{\ast \ast )}$ & $\left\vert \Gamma \right\vert =\sqrt{\det \Gamma }$                   \\ \hline
  norm & $\left\Vert q\right\Vert =\sqrt{a^{2}+b^{2}+c^{2}+d^{2}}$                    &   &                   \\ \hline
  &                    &  vigor  & vig$\Gamma = \Gamma \Gamma^{\ast}$                  \\ \hline
\end{tabular}

$^{\ast )}$or Clifford conjugation \hspace{1cm} $^{\ast \ast )}$ it exists on the set of proper or singular paravectors only!

 \vspace{5pt} For readers familiar with multivector algebra (Geometric Algebra) the following explanations are valid. Scalar and trivector are called a complex scalar and are denoted $\alpha=a+id$, where $i=a_{123}$ is a unity trivector. Vector and bivector are called a real and imaginary vector respectively, and their sum is a complex vector.     
             
In the next table we compare our notation with the multivector notation used by D. Hestenes in the work \cite{Hestenes}.

\begin{tabular}{|l|l||l|l|} 
\hline
         \multicolumn{2}{|c||} {Multivectors} & 
 \multicolumn{2}{|c|}{Paravectors} \\               
            \hline \hline
  multivector & $q=\left[ q\right] _{0}+\left[ q\right] _{1}+\left[ q\right]_{2}+\left[ q\right] _{3}$ & paravector & $\Gamma =
					\begin{bmatrix}
					a+id \\ 
					\textbf{b}+i\textbf{c}
					\end{bmatrix}$ \\
	\hline				
 scalar & $a=\left[ q\right]_{0}$ & real scalar & $a$ \\
 \hline
 vector & $\underline{\textbf{b}}=\left[ q\right] _{1}$ &real vector & $\textbf{b}$ \\
 \hline
bivector & $\underline{\textbf{C}}=\left[ q\right]_{2}$ & 
imaginary vector & $i \textbf{c}$ \\
\hline
trivector & $D=de_{123}=\left[ q\right] _{3}$ & imaginary scalar & $id$ \\
\hline
conjugation & $\overline{q}$ & reversion$^{\ast )}$ & $\Gamma ^{-}$ \\
\hline
reversion & $q^{\ast }$ & conjugation $^{\ast )}$ & $\Gamma ^{\ast }$ \\
\hline
involution & $\widehat{q}$ & & $\left( \Gamma^{-}\right) ^{\ast }$ \\
\hline
bivector coordinates & $C^{23}\mathbf{e}_{23},C^{31}\mathbf{e}_{31},C^{12}\mathbf{e}_{12}$ & & $ic^{1}\textbf{e}_{1},ic^{2} \textbf{e}_{2},ic^{3} \textbf{e}_{3},$ \\
\hline
trivector coordinates & $de_{123}$ & & $id$ \\
\hline
 & $\underline{\textbf{b}}_{1}\cdot \underline{\textbf{b}}_{2}$ & & $\textbf{b}_{1}\textbf{b}_{2}$ \\
 & $\underline{\textbf{b}}_{1}\cdot \underline{\textbf{C}}_{2}=-
\underline{\textbf{C}}_{2}\cdot \underline{\textbf{b}}_{1}$ & & $-\textbf{b}_{1}\times \textbf{c}_{2}$ \\
 inner products & $\underline{\textbf{C}}_{1}\cdot \underline{\textbf{C}}_{2}
                =\underline{\textbf{C}}_{2}\cdot \underline{\textbf{C}}_{1}$ & & $-\textbf{c}_{1}\textbf{c}_{2}$ \\
  & $D_{1}\cdot \underline{\textbf{b}}_{2}$  & & $id_{1}\underline{\textbf{b}}_{2}$ \\
  &  $D_{1}\cdot \underline{\textbf{C}}_{2}$ & & $-d_{1}\textbf{c}_{2}$ \\
  \hline
 \multirow{2}{*}{exterior products} & $\underline{\textbf{b}}_{1}\wedge \underline{\textbf{b}}_{2}$ & & $i\textbf{b}_{1}\times \textbf{b}_{2}$ \\
 & $ \underline{\textbf{b}}_{1}\wedge \underline{\textbf{C}}_{2}
                =\underline{\textbf{C}}_{2}\wedge \underline{\textbf{b}}_{1}$ & & $i\textbf{b}_{1}\textbf{c}_{2}$ \\
  \hline              
vector product & $\underline{\textbf{C}}_{1}\times \underline{\textbf{C}}_{2}
                =-\underline{\textbf{C}}_{2}\times \underline{\textbf{C}}_{1}$ & & $-i\textbf{c}_{1}\times \textbf{c}_{2}$   \\
                \hline                     
  
 \end{tabular}

		  *) Change of meaning between the conjugation and reversion was justified in section 4.

In matrix representation we can exhibit unit multivectors:

\qquad scalar\qquad $\qquad 
\begin{bmatrix}
1 &  &  & 0 \\ 
& 1 &  &  \\ 
&  & 1 &  \\ 
0 &  &  & 1
\end{bmatrix}
\qquad \qquad \qquad $vector\qquad $\begin{bmatrix}
0 & r_{x} & r_{y} & r_{z} \\ 
r_{x} & 0 & -ir_{z} & ir_{y} \\ 
r_{y} & ir_{z} & 0 & -ir_{x} \\ 
r_{z} & -ir_{y} & ir_{x} & 0
\end{bmatrix}
$

\qquad bivector$\qquad 
\begin{bmatrix}
0 & ir_{x} & ir_{y} & ir_{z} \\ 
ir_{x} & 0 & r_{z} & -r_{y} \\ 
ir_{y} & -r_{z} & 0 & r_{x} \\ 
ir_{z} & r_{y} & ir_{x} & 0
\end{bmatrix}
\qquad \qquad $trivector \ \ \ $\begin{bmatrix}
i &  &  & 0 \\ 
& i &  &  \\ 
&  & i &  \\ 
0 &  &  & i
\end{bmatrix}
$

where $\left\vert \mathbf{r}\right\vert =1.$

 \vspace{2pt}         
          \textbf{Paravectors in Pauli matrices basis}.           
                    
\vspace{2pt} Pauli matrices form the basis of complex 2x2 matrices in real space          
          
          $\sigma_{x}=\begin{bmatrix}0 & 1 \\ 1 & 0 \end{bmatrix} \qquad 
          \sigma_{y}=\begin{bmatrix}0 & -i \\ i & 0 \end{bmatrix} \qquad
         \sigma_{z}=\begin{bmatrix}1 & 0 \\ 0 & -1 \end{bmatrix} \qquad$
          
The basis in 4-dimensional space is obtained after adding matrix: $\sigma_{0}=\begin{bmatrix}1 & 0 \\ 0 & 1 \end{bmatrix}$
        
With the help of these matrices, we can show any paravector:

          $\alpha = (a+id)\sigma_{0}$

          $\pmb{\beta} = [(b_{x} +ic_{x})\sigma_{x} ,(b_{y} +ic_{y})\sigma_{y} , (b_{z} +ic_{z})\sigma_{z} ]$
        
         More information on the relationship between formalisms applied in physics of space-time can be found in the article \href{http://geometry.mrao.cam.ac.uk/1993/01/imaginary-numbers-are-not-real-the-geometric-algebra-of-spacetime/} {Imaginary Numbers are not Real.} The Geometric Algebra of Spacetime  Gull S., Lasenby A., Dorn Ch.
\section*{Acknowledgments}
I give great thanks to my daughter, Dorota, for her careful reading of this manuscript and her help in my mastering of TEX editor.            

\bibliographystyle{plain}
\bibliography{bibliografie}

\end{document}